\documentclass[a4paper]{article}
\usepackage{geometry}
\usepackage[utf8]{inputenc}
\usepackage{authblk}
\usepackage{physics}
\usepackage{ragged2e}
\usepackage{mathtools}
\usepackage{amsfonts} 
\usepackage{amsthm}
\usepackage[style=alphabetic]{biblatex}
\renewenvironment{abstract}{%
   \MakeUppercase{\abstractname.}}{%
  \vspace*{\baselineskip}}
  
\newtheorem{teo}{Theorem}[section]
\newtheorem{lem}[teo]{Lemma}
\newtheorem{prop}[teo]{Proposition}

\newtheorem{con}{Conjecture}

\theoremstyle{definition}
\newtheorem{rem}[teo]{Remark}

\def\C{\mathbb{C}}
\def\R{\mathbb{R}}
\def\Q{\mathbb{Q}}
\def\Z{\mathbb{Z}}
\def\N{\mathbb{N}}

\def\D{\Delta}
\def\h{\hat{h}}

\title{\textbf{Counting rational points on elliptic curves with a rational 2-torsion point}}
\author{FRANCESCO NACCARATO\thanks{Faculty of Sciences, Scuola Normale Superiore, Pisa. Contact: francesco.naccarato@sns.it.}}

\date{}

\begin{document}

\maketitle

\begin{abstract}
Let $E/\Q$ be an elliptic curve over the rational numbers. It is known, by the work of Bombieri and Zannier, that if $E$ has full rational $2$-torsion, the number $N_E(B)$ of rational points with Weil height bounded by $B$ is $\exp (O\left(\frac{\log B}{\sqrt{\log\log B}}\right))$. In this paper we exploit the method of descent via $2$-isogeny to extend this result to elliptic curves with just one nontrivial rational $2$-torsion point. Moreover, we make use of a result of Petsche to derive the stronger upper bound $N_{E}(B) = \exp(O\left(\frac{\log B}{\log\log B}\right))$ for these curves and to remove a deep transcendence theory ingredient from the proof.
\end{abstract}

\section{Introduction and results}\par

For an elliptic curve $E/\Q$ in Weierstrass form given by an affine equation:
\begin{align} y^2 = x^3+ax^2+bx+c \ \ \ \ \ \ a,b,c \in \Q \label{eq:1}\end{align}
we define its (naive) height $H(E)$ as the Weil height of the vector $(1,a,b,c)$, its discriminant $\D_E$ as $16$ times the discriminant of the polynomial on the right-hand side (which is nonzero by the nonsingularity of elliptic curves) and we let $h(E):=\log H(E)$.\par

We recall that the points of an elliptic curve defined over $\C$ form a group with the usual structure arising from the Weierstrass map. When $E$ is defined over $\Q$, the set of its rational points is a finitely generated subgroup, the Mordell-Weil group $E(\Q)$; the rank $r_E$ is defined as the abelian rank of $E(\Q)\simeq\Z^{r_E}\times T$, where $T$ is the torsion subgroup.\par
Let us introduce the usual Weil height $H(P)$ of a rational point as the Weil height of its $x$-coordinate, and  its logarithmic analogue as $h(P) := \log H(P)$. Moreover, let us define the quantity that we are interested in bounding:
\begin{align*} N_E(B):=|\{P \in E(\Q):H(P) \le B\}|\end{align*}
or, equivalently, \begin{align} N_E(B):=|\{P \in E(\Q):h(P) \le \log B\}|.\label{eq:4}\end{align}\par
We will from now on omit the subscript $E$ on $\D$, $r$ and $N(B)$, keeping in mind the dependence of these quantities on the curve.\par
\vspace{3mm}
In \cite{izvestiya}, Theorem 1, Bombieri and Zannier proved that for elliptic curves in Weierstrass form with full rational $2$-torsion one has \begin{align*}N(B) \le B^{\frac{c}{\sqrt{\log\log B}}}\end{align*} for sufficiently large $B$ (depending on the curve), with $c$ an absolute constant.\par\vspace{3mm}
We can now state our result:
\begin{teo}There exists absolute computable constants $C, \ c_0$ such that for any elliptic curve 
$E/\Q$ as in \eqref{eq:1} with a rational $2$-torsion point the inequality \begin{align}N(B) \le 
B^{\frac{C}{\log\log B}}\label{eq:5}\end{align} holds for all $B\ge \max\{e^e, (eH(E))^{c_0}\}$.\end{teo}\par
This result implies, for this special family of curves, a conjecture that is widely believed to hold true for any elliptic curve over $\Q$ in Weierstrass form:
\begin{con}Let $\epsilon > 0$. There exists a constant $c'(\epsilon)$ depending only on 
$\epsilon$ such that \begin{align*}N(B) \le 
c'(\epsilon)B^{\epsilon}\end{align*} for any elliptic curve $E/\Q$ as in \eqref{eq:1} for sufficiently large $B$ again depending on the equation of the curve.\end{con}
We are able to obtain our improved upper bound for $N(B)$ by making use of a result of Petsche \cite{petsche}, which we will state in Section 4, giving a lower bound for the smallest height of a non-torsion point that depends polynomially on the \textit{Szpiro ratio} of the curve.
By contrast, the lower bound for the same quantity used in \cite{izvestiya} depends exponentially on the \textit{Szpiro ratio} (see Theorem 0.3 in Hindry and Silverman's work \cite{HS}), and turns out to be useful only if the rank of the curve is large enough. For small ranks, Bombieri and Zannier employed a result of Masser (see \cite{masser}, Theorem), whose proof relies on transcendence theory techniques. Hence, by removing the need for this result, we are considerably simplifying the proof structure.\par
Furthermore, as we will see in Section 5, even assuming Lang's height conjecture for rational elliptic curves, which improves on Petsche's estimate for the smallest canonical height by removing the dependence on the \textit{Szpiro ratio}, does not lead to an improvement on our estimate.
\par

The need for a single $2$-torsion point comes from the fact that it turns out to be sufficient
to effectively bound the rank of the curve in terms of the discriminant, as we will illustrate
in Section 3.

\section{Canonical height and quasi-minimal model}
\par
Recall that a Weierstrass equation for $E$ is called a \textit{minimal model} at a prime $p$ if the $p$-adic valuation of the relative discriminant is the smallest possible (that is, an integer $0 \le m <12$). An equation which is a minimal model at all primes is called a (global) minimal model.\par
We will see that it is not restrictive to work with a quasi-minimal model (one that is minimal at all primes $p$ except at most $2$ and $3$) of $E$, given by:
\begin{align} y^2=x^3+Ax+B \label{eq:2}\end{align}
with $A,B \in \Z$. We see that in this case $\D = -16(4A^3+27B^2)$, so we have: \begin{align}|\D| \le 496H(E)^3.\label{eq:3}\end{align}

For our estimation, we make use of the canonical height:
\begin{align*}\h(P):=\lim_{n \rightarrow \infty}4^{-n}h(2^nP).\end{align*}
Set \begin{align*}N_{can}(B):=|\{P \in E(\Q):\h(P) \le \log B\}|.\end{align*}
We will derive Theorem 1.1 from the following:
\begin{teo}There exists an absolute computable constant $D$ such that for any elliptic curve $E/\Q$ as in \eqref{eq:2} with a rational $2$-torsion point the inequality \begin{align}N_{can}(B) \le 
B^{\frac{D}{\log\log B}}\label{eq:6}\end{align} holds for all $B \ge \max\{e^e, H(E)\}$.\end{teo}

In what follows, the numbers $c_k$ will be positive, absolute and computable constants. We will prove Theorem 2.1 in the next two sections; let us first derive our main result from it. We have:

\begin{prop}Theorem 2.1 implies Theorem 1.1.\end{prop}
\textbf{Proof}: First, notice that inequality $\eqref{eq:6}$ for elliptic curves as in $\eqref{eq:2}$ extends to elliptic curves as in $\eqref{eq:1}$, with a different lower bound for $B$: given $E$ as in $\eqref{eq:1}$, the usual translation-dilation isomorphism $\phi: E \rightarrow E'$ to the respective model as in \eqref{eq:2} preserves the canonical height and we have $h(E') \le c_0h(E)+c_0$ (this is the $c_0$ in the statement of Theorem 1.1; let us also take $c_0\ge1$), so for $E$ we obtain $\eqref{eq:6}$ for all $B\ge\max\{e^e, (eH(E))^{c_0}\}$.\par
On the other hand, Zimmer \cite{zimmer} proved that, for a fixed curve in Weierstrass model\\ $Y^2=4X^3-g_2X-g_3$, the difference between the Weil and canonical 
heights is small:\\
$|\h(P)-h(P)|<c_1h(E)+c_1$. As remarked in \cite{izvestiya}, this is easily extended to models as in \eqref{eq:1}. Since by hypothesis $\log\log B \ge 1$ and $B\ge (eH(E))^{c_0}\ge H(E)$, we have
\begin{align*}N(B)\le N_{can}(\exp(\log B+c_1h(E)+c_1)) \le (B(eH(E))^{c_1})^{\frac{D}{\log\log B}}\le B^{\frac{(2c_{1}+1)D}{\log\log B}}.\end{align*} This completes the proof.
\vspace{3mm}

It is well known that $\h$ is a quadratic form on $E(\Q)$ and that it is positive definite on the quotient $E(\Q)/T \simeq \Z^r$. If we consider the tensor product of abelian groups $E(\Q) \otimes \R \simeq \R^r$, we see that the torsion is mapped in $0$ by the tensor, and $E(\Q)/T$ injects (in the sense of $E(\Q)/T \otimes 1 \hookrightarrow E(\Q) \otimes \R$) in a lattice $\mathcal L_E$ of dimension $r$: we can then bring our quadratic form onto this lattice by setting $\h(x):=\h(P)$ where $P$ is the unique rational point modulo torsion such that $P\otimes1=x$.\par
This form can be extended by linearity to $\Q^r$ and by continuity to all $\R^r$, and it is a well-known result that it remains positive definite: hence, we can take $\h$ to be the square of the usual euclidean norm (so $\mathcal L_E$ won't necessarily be $\Z^r$ in a metric sense).\par
The problem of counting rational points with Weil height up to $B$ is thus equivalent to that of counting the number of points of $\mathcal L_E$ in the ball of centre $0$ and radius $\sqrt{\log B}$, and then multiplying by the cardinality of the torsion. We immediately remark that for elliptic curves over $\Q$ the cardinality of the torsion is known to be absolutely bounded by the work of Mazur \cite{mazur}.\par
In the light of the above, from now on we will work with the canonical height, omitting the subscript from $N_{can}(B)$.

\section{Bounding the rank}
We now turn to the proof of Theorem 2.1. As we will see in the next section, in order for our methods to work on a given elliptic curve it is crucial to have an upper bound on its rank. If our curve has a rational $2$-torsion point, a good upper bound in terms of the discriminant can be obtained in a simple way through a descent via $2$-isogeny.
\begin{lem}Let $E/\Q$ be an elliptic curve as in \eqref{eq:2} with a rational $2$-torsion point. Let $\D$ be its discriminant and $r$ its rank: then one has \begin{align*}r \le 2\omega(\D)+2\end{align*}
where $\omega(m)$ is the number of distinct prime factors of the integer $m$.\end{lem}
\textbf{Proof}: We can apply a descent via $2$-isogeny as described in \cite{silverman1} (see specifically p.92, Proposition 3.8 and p.98), obtaining: \begin{align*}|E(\Q)/2E(\Q)| \le 2^{\omega(\D)+2}.\end{align*}
This, together with the bound \begin{align*}2^r \le |E(\Q)/2E(\Q)|,\end{align*} which follows by explicitly writing the quotient as \begin{align*}\left(\Z/2\Z\right)^r\prod_{i=1}^{k}\Z_{p_i^{r_i}}/2\Z_{p_i^{r_i}}\end{align*} thanks to the classification of finitely generated abelian groups, gives:
\begin{align*}r \le 2\omega(\D)+2. \ \ \square\end{align*}

Observe that \eqref{eq:3} tells us that \begin{align}\log B \ge \max\{e, h(E)\} \ge \max\{e, \frac{1}{3}\log |\D|-\log 8\} \ge \frac{1}{6}\log |\D|.\label{eq:7}\end{align} 
Since the \textit{Prime Number Theorem} implies, for $m\ge3$, the bound \begin{align*} \omega(m) \le c_2\frac{\log |m|}{\log\log |m|}\end{align*}
and since $|\D|\ge16\ge3$, from Lemma 3.1 and \eqref{eq:7} we obtain: \begin{align}r\le c_3\frac{\log B}{\log\log B}.\label{eq:8}\end{align}

Let us remark that the existence of a rational $2$-torsion point is used exclusively in order to obtain Lemma 3.1.

\section{Counting through a covering argument in $\R^n$}

We will apply the following well-known counting strategy:

\begin{enumerate} \item we find a small enough radius $\rho_0$ such that we can ensure that any ball of radius $\rho_0$ centered at a point of $\mathcal L_E$ does not contain any other point of $\mathcal L_E$;
\item we count how many of these balls we need to cover the intersection of $\mathcal L_E$ with the ball of centre $0$ and radius $\sqrt{\log B}$. \end{enumerate}

This is almost the same strategy followed by Bombieri and Zannier in \cite{izvestiya}: the only substantial difference is that they work with two different values of $\rho_0$ depending on the magnitude of the rank, and for one of these values the number of lattice points inside the small balls is not necessarily $1$, but some computable constant.\par
For the reader's convenience, we repeat the lemmas and proofs that can already be found in \cite{izvestiya}.
\vspace{3mm}

The second step of the strategy is the easiest, for we just need an elementary covering lemma:
\begin{lem}Given a positive integer $n$, radii $R, \rho$ and a subset $S$ of the $n$-ball $\mathcal B_n(0,R)$, there exists a set of at most $(1+\frac{2R}{\rho})^n$ balls of radius $\rho$ centered at points of $S$ such that $S$ is contained in the union of these balls.\end{lem}
\textbf{Proof}: Consider a maximal set $\Gamma$ of disjoint balls of radii $\frac{\rho}{2}$ centered at points of $S$ (by maximal we mean such that any ball of radius $\frac{\rho}{2}$ centered at a point of $S$ intersects a ball of $\Gamma$). Notice that the union of the balls in $\Gamma$ is contained in $\mathcal{B}_n(0, R+\frac{\rho}{2})$, so we necessarily have: \begin{align*}|\Gamma| \le \frac{V_n(R+\frac{\rho}{2})}{V_n(\frac{\rho}{2})} \le \left(1+\frac{2R}{\rho}\right)^n\end{align*} with $V_n(a)$ the $n$-volume of the $n$-ball of radius $a$.\par
By enlarging the balls in $\Gamma$ by a factor of $2$ we get a set of balls of radius $\rho$ centered at points of $S$ covering $S$, because if a point of $S$ lied outside the union of these new balls, that would contradict maximality. $\square$
\vspace{3mm}

It is clear that we will make use of this lemma putting $n=r$ (and $R=\sqrt{\log B}$, $\rho=\rho_0$), from which the importance of the magnitude of the rank in our estimation follows.\par

A good value for $\rho_0$ is derived by Petsche in \cite{petsche}. It turns out that the minimum canonical height of a non-torsion point can be bounded below as a function of the minimal discriminant $\mathcal D$ of $E$ (that is, the discriminant of one of its minimal models) and of its \textit{Szpiro ratio}:
\begin{align*}\sigma=\frac{\log |\mathcal D|}{\log{\mathcal N}}\end{align*}
where $\mathcal N$ is the conductor of the curve.
\begin{rem}If $\mathcal{N} = 1$ then it is known that $\mathcal D=1$ and $\sigma$ is defined to be $1$.\end{rem}

\begin{rem}Recall that the prime factors of the minimal discriminant of an elliptic curve $E/\Q$ are exactly the primes of bad reduction for $E$, whereas the discriminant $\D$ of our quasi-minimal model has an additional factor $2^j3^k$ with $j,k \in \N$. For the definition of the conductor and the proof that $\sigma \ge 1$ see, for example, \cite{petsche}. \end{rem}

\begin{prop}[Petsche, \cite{petsche}] There exist constants $c_4, c_5 >1$ such that for any non-torsion point $P \in E(\Q)$, we have: \begin{align*}\h(P) \ge \frac{\log |\mathcal D|}{c_4\sigma^6\log^2(c_5\sigma)}.\end{align*} \end{prop}

In Section 5 we will see that this can be slightly strengthened for rational elliptic curves, but the improvement is not relevant for our estimates.\par
Notice that $\log^2(c_5\sigma) < c_5^2\sigma^2$, so we have: 
\begin{align}\h(P) > \frac{\log |\mathcal D|}{c_6\sigma^8}.\end{align}\par

Observe that since no elliptic curve over $\Q$ has everywhere good reduction, $|\mathcal{D}|>1$ and our lower bound is positive. For the same reason, Remark 4.2 implies that $\mathcal{N} \ge 2$. Let us now apply Lemma 4.1 with $\rho=\sqrt{\frac{\log |\mathcal D|}{c_6\sigma^8}}$; then we know that in each of the balls of this radius centered at points of $\mathcal L_E$ there is just one non-torsion point (that is, the centre) and so we obtain: \begin{align} N(B) \le |T|\left(1+2c_6\sqrt{\frac{\log B}{\log |\mathcal D|}}\sigma^4\right)^r \le c_7\left(c_8\frac{\log B}{log\mathcal N}\sigma^3\right)^{r}\label{eq:11}\end{align}
since surely $|\mathcal D| \le |\D|$. \par
Observe that for the set of curves of rank bounded by any fixed absolute constant $c_9$ the statement of Theorem 2.1 clearly follows from \eqref{eq:11}, since in this case we can just directly bound $\log \mathcal{N} \ge \log 2$, $\sigma \le 2\log |\mathcal D| \le 2c_2\log B$ and get:
\begin{align*} N(B) \le c_7(c_{10}\log B)^{4c_9} \le (\log B)^{c_{11}}\le B^{\frac{c_{11}}{\log \log B}}.\end{align*}
So from now on we can consider, for example, $r \ge 39$. Notice that this in turn requires $B$ to be not too small.\par
Since, by \cite{silverman2}, Appendix C, Table 15.1, $\mathcal N$ is divisible by all the primes dividing $\mathcal D$ except at most $2$ and $3$, and hence by all the primes dividing $\Delta$ except at most $2$ and $3$, it follows from Lemma 3.1 that $6\mathcal N$ has at least $\max\{2,\frac{r-2}{2}\} \ge \frac{r}{3}$ prime factors and hence we have:
\begin{align}6\mathcal{N} \ge p_{\frac{r}{3}}\# \label{eq:12}\end{align} with $p_n\# := \prod_{j=1}^{n}p_j$ the product of the first $n$ primes.

\begin{lem}For $n \ge 13$ the inequality: \begin{align*} p_n\# \ge n^n\end{align*} holds. \end{lem}

\textbf{Proof}: This is easily verified with a calculator for $n=13$, and then we can proceed by induction: for the inductive step we just need $p_n \ge \frac{n^n}{(n-1)^{n-1}}$. But $\frac{n^n}{(n-1)^{n-1}} \le n\left(1+\frac{1}{n-1}\right)^{n-1} \le en$ and $p_n > n\log n > en$ by Rosser's theorem \cite{ross} and the fact that $p_n \ge 43 > e^e$. $\ \square$
\vspace{3mm}

Lemma 4.5 and \eqref{eq:12} give us: 
\begin{align*}\log 6 + \log\mathcal{N} \ge \frac{r}{3}\log \frac{r}{3}.\end{align*}
For $r \ge 39$ we obtain $\log\mathcal{N} \ge \frac{r\log r}{7}$ and consequently $\sigma \le \frac{7c_2\log B}{r\log r}$. Substituting in \eqref{eq:11} we are left with the task of showing that the maximum of the function 
\begin{align*}f(x)=c_7\left(c_{12}\frac{\log B}{x\log x}\right)^{c_{13}x}\end{align*} for $x \in [39, \ c_3\frac{\log B}{\log \log B}]$ is $B^{O\left(\frac{1}{\log \log B}\right)}$. This is done as follows: clearly we can forget about the constants $c_7$ and $c_{13}$. We can also forget about $c_{12}$ precisely because $x = O\left(\frac{\log B}{\log\log B}\right)$. \par
Write $A=\log B$, $f(x)=\exp(x(\log A - \log(x\log x)))$.
The derivative \begin{align*}f'(x)=f(x)[\log A - (\log(x\log x) + 1 + \frac{1}{\log x})]=:f(x)g(x)\end{align*} is positive for $x$ up to a certain value $x_0$ that satisfies $x_0\log x_0 \sim \frac{A}{e}$ and negative afterwards (since $x \ge 39$ we don't have to worry about the term $\frac{1}{\log x}$ being big). Formally, for $x \ge 39$, $g(x)$ is a decreasing function in $x$, which is nonnegative for $x\log x=Ae^{-1-\frac{1}{\log 39}}$ and negative for $x\log x=A$.
This in particular tells us that the zero $x_0$ of $f'(x)$ satisfies $x_0 \le 2\frac{A}{\log A}$ and hence: \begin{align*}f(x) \le \exp(2\left(1+\frac{1}{\log 39}\right)\frac{A}{\log A}) = B^{O\left(\frac{1}{\log \log B}\right)}\end{align*}
as required.

\section{Some remarks and further explorations}

Observe that, even if we had a lower bound
\begin{align*}\h(P) \gg \log |\mathcal D|,\end{align*} independent of $\sigma$, which is Lang's height conjecture for elliptic curves over $\Q$ (or even just for those with a rational $2$-torsion point), we could not improve on Theorem 2.1 by means of our methods: in fact, we would have to bound \begin{align*}\left(\frac{\log B}{\log |\mathcal D|}\right)^r\end{align*} for $\log |\mathcal D|$ and $r$ in the ranges given by \eqref{eq:3} and \eqref{eq:8}. Choosing $\log B = c_{14}\log |\mathcal D|$ and clearly $r$ as big as possible, we get exactly the result we already proved:
\begin{align*}N(B) \le e^{c_{15}\frac{\log B}{\log\log B}} = B^{\frac{c_{15}}{\log\log B}}.\end{align*}\par
This hints that, in order to improve on our bound for $N(B)$, one should improve on the estimate for the rank in terms of the discriminant (at least for our special class of curves; moreover, obtaining any sensible estimate of this kind for a broader class of curves could help in establishing Conjecture 1 for more elliptic curves). As an extreme example, observe that uniform boundedness of the ranks would give an upper bound of $(\log B)^{O(1)}$ for $N(B)$.
\vspace{4mm}

We point out that our upper bound for $N(B)$ is barely insufficient for what would have been an interesting application: in \cite{BB} Bombieri and Bourgain look at the diophantine equation $x^2+y^2=m$. Specifically, they are interested in finding an asymptotic upper bound to the cardinality of the set \begin{align*} S_6(m) = \{(\lambda_1,...,\lambda_6) \in \Lambda_m^6: \lambda_1+\lambda_2+\lambda_3=\lambda_4+\lambda_5+\lambda_6\}\end{align*}
where $\Lambda_m$ is the set of the gaussian integers $\lambda$ such that $\lambda\bar{\lambda}=m$. Setting $N = |\Lambda_m|$ the number of solution to the equation of interest, their goal is to prove that one has: $|S_6(m)| = O(N^{3+\epsilon}) \ \forall \epsilon > 0$.\par
Tackling this problem through the theory of elliptic curves, they are able to prove a weak conditional result (Theorem 8), that we restate here for the reader's convenience:
\begin{teo}[Bombieri, Bourgain] At least one of the two following statements holds.
\begin{enumerate}
    \item If $\frac{\log N}{\log\log m} \rightarrow \infty$, then $|S_6(m)|=O(N^{3+o(1)})$.
    \item There exist elliptic curves $E/\Q$ of unbounded rank.
\end{enumerate}\end{teo}
It would not be free of interest to remove the condition on the ranks, which appears to be far from easy to settle. Since the curves they deal with have a rational $2$-torsion point, we can apply our result to their methods and remove the dependence from the rank. By doing this, one sees that the first condition has to be weakened to $\frac{\log N\log\log m}{\log m} \rightarrow \infty$, which turns out to be vacuous because $\log N = O\left(\frac{\log m}{\log\log m}\right)$, as can be inferred by writing explicitly $N$ as in \cite{BB}, (20).\par
Any asymptotic improvement on our bound would instead imply the conjecture for 
the family: \begin{align*}\{m : \log N \ge A\frac{\log m}{\log\log 
m}\}\end{align*} for any fixed $A>0$, which for sufficiently small $A$ is 
nonempty (it contains the products of the first $k$ primes 
congruent to $1$ mod $4$, for example).
\vspace{4mm}

Even though it does not lead to any quantitative improvement in our estimate, other than the change of an absolute constant, we remark that Petsche's lower bound for the smallest canonical height (which in the original paper is established for curves over arbitrary number fields and depends also on the degree of the number field) can be slightly improved for curves over $\Q$. More precisely, in this case one has: \begin{align*}\h(P) \ge \frac{\log |\mathcal D|}{c_{16}\sigma^6}.\end{align*}\par
To see this, consider Proposition 8 of \cite{petsche}, that we recall for the reader's convenience:
\vspace{2mm}

Let $k$ a number field of degree $d = [k : \Q]$, and let $E/k$ be an elliptic curve with Szpiro ratio $\sigma$. Then: \begin{align*}|\{P \in E(k) : \h(P) \le \frac{\log{\N_{k/\Q}(\D_{E/k})}}{2^{13}3d\sigma^2}\} \le a_1d\sigma^2\log{(a_2d\sigma^2)}\end{align*}
with $a_1=134861$ and $a_2=104613$.
\vspace{2mm}

In its proof in the original paper, we see that in the precise case of $k=\Q$ that we are examining, the left-hand side of (25) is precisely $0$ since $\Q$ has just one archimedean place, and no estimation is needed. This enables us to remove the $\log{(a_2d\sigma^2)}$ factor in the statement of Proposition 8. The rest of the proof remains unchanged.\par
As remarked by Petsche, this also holds for imaginary quadratic fields, since they, too, have just one archimedean place.

\vspace{6mm}

\textbf{\large{Acknowledgments}}
\vspace{2mm}

The author would like to thank Professor Umberto Zannier for the invaluable
guidance and suggestions, not only of a mathematical nature, which he 
offered throughout the process of writing this paper.\par
The author would like to thank Professor Clayton Petsche for the 
helpful discussion concerning his work.\par
The author would also like to thank the anonymous referee for the careful and insightful review of the paper.

\printbibliography
\end{document}